\documentclass[12pt]{article}

\usepackage{t1enc,amsmath,amssymb,amsthm,amscd,epsfig,latexsym,amsfonts,ifthen}

\usepackage{amsmath}%
\usepackage{amsfonts}%
\usepackage{amssymb}%
\usepackage{graphicx}
\newtheorem{theorem}{Theorem}
\theoremstyle{plain}

\newtheorem{corollary}{Corollary}

\newtheorem{lemma}{Lemma}

\numberwithin{equation}{section}

\begin{document}

\title{\vspace{-1in}\parbox{\linewidth}{\footnotesize\noindent}
 \vspace{\bigskipamount} \\
The Cauchy-Pompeiu integral formula in elliptic complex numbers
\thanks{ 
\hfil\break\indent
{\em $\dagger$ danieldaniel@gmail.com ~
\em $\ddagger$ cvanegas@usb.ve }}}

\author{D. Alay\'on-Solarz$\dagger$ and C. J. Vanegas$\ddagger$ \\
Departamento de Matem\'aticas Puras y Aplicadas \\
Universidad Sim\'{o}n Bol\'{\i}var, Caracas 1080-A, Venezuela}

\date{}

\maketitle

\begin{abstract}

The aim of this article is to give a generalization of the Cauchy-Pompeiu integral 
formula for functions valued in pa\-ra\-me\-ter-depending elliptic algebras with  structure 
polynomial $X^2 + \beta X + \alpha$ where $\alpha$ and $\beta$ are real numbers. As a consequence, a Cauchy integral representation 
formula is obtained for a generalized class of holomorphic functions.

\end{abstract}

\section{Introduction}

In his book \cite{Yag}, I. Yaglom introduced the "Generalized Complex Numbers", as complex numbers of the form $z = x + iy$ where the product of two complex numbers is induced by the relation
\begin{displaymath}
i^2 = - \beta i - \alpha
\end{displaymath}
where $\alpha$ and $\beta$ are real numbers  subjected to the ellipticity condition  $4 \alpha - \beta^2 > 0$. Yaglom observed 
that generalized complex numbers are isomorphic to the ordinary complex numbers and thus from the point of view 
of the algebra they are structurally no different from the latter. 

In the language developed in \cite{Tut2} a general complex algebra in the sense of Yaglom is a 
parameter-dependent $2$-dimensional algebra type Clifford  with structure polynomial
\begin{displaymath}
X^2 + \beta X + \alpha
\end{displaymath}
and subjected to the condition that $\alpha$ and $\beta$ are real constants satisfying $4 \alpha - \beta^2 > 0$. 

In contrast to the strictly algebraic point of view, in Analysis a significant gain is obtained when we consider the Cauchy-Riemann operator
\begin{displaymath}
\partial_{\bar z} = \frac{1}{2}(\frac{\partial }{\partial x} + i \frac{\partial}{\partial y})
\end{displaymath}
acting on a complex valued function $f(z) = u + i v$ in the generalized complex algebra determined by $\alpha$ and $\beta$, and considering functions in the kernel of this operator as the holomorphic ones. As there exist differentiable functions that are not holomorphic in the ordinary sense, yet they are holomorphic for some suitable choice of real numbers $\alpha$ and $\beta$, a more encompassing concept of holomorphicity is obtained.

The purpose of this article is to enunciate and prove a generalization of the ordinary Cauchy-Pompeiu integral formula for complex algebras with structure polynomial $X^2 + \beta X + \alpha$ of elliptic type and where $\alpha$ and $\beta$ are constants, given by
\begin{displaymath}
f(\zeta) = \frac{1}{ 2 \pi \hat{i} }  \int\limits_{\partial \Omega} \frac{f(z)}{\widetilde{z - \zeta}} d \tilde z - \frac{1}{ \pi \hat{i}}  \iint\limits_{ \Omega} \frac{\partial_{\bar z} f(z)}{\widetilde{z - \zeta}} dx dy 
\end{displaymath}
where
\begin{displaymath}
\ \ \ \  \tilde{z} = y - ix , \ \ \ d \tilde{z} = dy - i dx\ \ \  \hat{i} = \frac{\beta + 2 i}{\sqrt{4 \alpha - \beta^2 }}.
\end{displaymath}
In particular, if $\alpha = 1$ and $\beta = 0$ we have that
\begin{displaymath}
\tilde z = \frac{z}{i} , \ \ \  d \tilde z = \frac{1}{i} dz, \ \ \ \ \hat{i} = i
\end{displaymath}
and we obtain
\begin{displaymath}
f(\zeta) = \frac{1}{ 2 \pi i }  \int\limits_{\partial \Omega} \frac{f(z)}{z - \zeta} d z - \frac{1}{ \pi }  \iint\limits_{ \Omega} \frac{\partial_{\bar z} f(z)}{z - \zeta} dx dy,
\end{displaymath}
recovering the ordinary Cauchy-Pompeiu representation formula. An immediate consequence is the Cauchy integral formula for holomorphic functions in the generalized complex algebra given by
\begin{displaymath}
f(\zeta) = \frac{1}{ 2 \pi \hat{i} }  \int\limits_{\partial \Omega} \frac{f(z)}{\widetilde{z - \zeta}} d \tilde z.
\end{displaymath}
\section{Preliminaries}

Let $z = x + i y$ denote a complex number in the Clifford Algebra with structure polynomial $X^2 + \beta X + \alpha$, where $\alpha$ and $\beta$ are real numbers. In this algebra the product is defined for two complex numbers $z_1 = x_1 + i y_1$ and $z_2 = x_2 + i y_2$ as
\begin{displaymath}
z_1 z_2 = (x_1 x_2 - \alpha y_1 y_2) + i (x_1 y_2 +y_1 x_2 - \beta y_1 y_2).
\end{displaymath}

It can be shown that the product is uniquely invertible provided that $4 \alpha - \beta^2 > 0$, as an inverse is given by
\begin{displaymath}
z^{-1} = \frac{x - \beta y - i y}{x^2 - \beta x y + \alpha y^2 }
\end{displaymath} 
is well defined for every $z \neq 0$. We call this type of algebra to be of elliptic type. In this case
the expression
\begin{displaymath}
<z , w>_{(\alpha, \beta)} = \frac{1}{2}( w \cdot ( \bar{z} - \beta \cdot \text{Im} (z)) +  ( \bar{w} - \beta \cdot \text{Im} (w)) z)
\end{displaymath}
where $\bar z = x - i y$ is the conjugation and  $\text{Im}( w)$ and $\text{Im}( z)$ denote the imaginary part of $w$ and $z$ respectively, defines an inner product over the real numbers. Then a norm can be defined as
\begin{displaymath}
||z||_{(\alpha, \beta)} = \sqrt{<z,z>_{(\alpha, \beta)}} = \sqrt{x^2 - \beta x y + \alpha y^2}.
\end{displaymath}
In particular $|| \cdot ||_{(0, 1)}$ is the Euclidean norm. The distance depending on the numbers $\alpha$ and $\beta$ is compatible with the algebra with structure polynomial $X^2 + \beta X + \alpha$ in the following sense, for all $w$ and $z$ complex numbers we have
\begin{displaymath}
|| z \cdot w ||_{(\alpha, \beta)}  = ||z||_{(\alpha, \beta)} ||w||_{(\alpha, \beta)} ,
\end{displaymath} 
and for every $z \neq 0$
\begin{displaymath}
|| z^{-1}  ||_{(\alpha, \beta)} = \frac{1}{|| z  ||_{(\alpha, \beta)}}.
\end{displaymath}

The Cauchy-Riemann operator $\partial_{\bar z}$ is defined as
\begin{displaymath}
\partial_{\bar z} = \frac{1}{2}( \frac{\partial}{ \partial x} + i \frac{\partial }{\partial y})
\end{displaymath}
and the expression $\partial_{\bar z} w = 0$ where $w = u + i v$ is equivalent to the following real system of equations
\begin{displaymath}
\frac{1}{2}(\partial_x u - \alpha \partial_y v) = 0, ~~~~  \frac{1}{2}( \partial_y u +  \partial_x v - \beta \partial_y v) = 0,
\end{displaymath}
which is a generalization of the ordinary Cauchy-Riemann equations. In this case we say the function $w$ is holomorphic with respect to the structure polynomial $X^2 + \beta X + \alpha$.
 
\section{Five preliminary lemmas}

In order to be able to adapt the classical proof of the Cauchy-Pompeiu representation formula when $\alpha =1$ and $\beta = 0$ 
\cite{tut} to the case when $\alpha$ and $\beta$ are constants satisfying $4 \alpha - \beta^2 >0$, 
the following elementary lemmas will be useful:
\begin{lemma}
In every complex algebra with structure polynomial $X^2 + \beta X + \alpha$ with  $\alpha$ and $\beta$ constants, the product rule holds for the Cauchy-Riemann operator.
\end{lemma}
\proof Direct computation show that for every pair of differentiable functions $f_1$ and $f_2$ it holds
\begin{displaymath}
\partial_{\bar z} (f_1 \cdot f_2) = \partial_{\bar z} f_1 \cdot f_2 + f_1 \cdot \partial_{\bar z} f_2,
\end{displaymath}
 where the product is understood to be in the complex algebra with structure polynomial $X^2 + \beta X + \alpha $ with $\alpha$ and $\beta$ constants. 
 \qedhere
\begin{lemma}
For every algebra with structure polynomial $X^2 + \beta X + \alpha$ such that $4 \alpha - \beta^2 >0$ and $\alpha$ and $\beta$ are constants, the map $\big( \widetilde{z- \zeta} \big) ^{-1}$ is weakly singular at $\zeta$ in the ordinary sense.
\end{lemma}
\proof From the equivalence of the norms $|| \cdot ||_{(\alpha, \beta)}$ and $|| \cdot ||_{(1, 0)}$, there exist positive numbers $K_1(\alpha, \beta)$ and  $K_2(\alpha, \beta)$  such that for a fixed $\zeta$ and every $ z \neq \zeta$:
\begin{displaymath}
K_1(\alpha, \beta) || \big( \widetilde{z- \zeta} \big) ^{-1} ||_{(\alpha, \beta)} \leq || \big( \widetilde{z- \zeta} \big) ^{-1} ||_{(1, 0)} \leq  K_2(\alpha, \beta) || \big( \widetilde{z- \zeta} \big) ^{-1} ||_{(\alpha, \beta)}
\end{displaymath}
and 
\begin{displaymath}
K_1(\alpha, \beta) || \big( \widetilde{z- \zeta} \big) ||_{(\alpha, \beta)} \leq || \big( \widetilde{z- \zeta} \big)  ||_{(1, 0)} \leq  K_2(\alpha, \beta) || \big( \widetilde{z- \zeta} \big)  ||_{(\alpha, \beta)}.
\end{displaymath}
From the second line we infer that
\begin{displaymath}
\frac{1}{ || \big( \widetilde{z- \zeta} \big) ||_{(\alpha, \beta)} }\leq   K_2(\alpha, \beta) \frac{1}{ || \big( \widetilde{z- \zeta} \big)  ||_{(1,0)}} =  K_2(\alpha, \beta) \frac{1}{ || \big( z- \zeta \big)  ||_{(1,0)}},
\end{displaymath}
and since 
\begin{displaymath}\frac{1}{|| \big( \widetilde{z- \zeta} \big) ||_{(\alpha, \beta)}} =  || \big( \widetilde{z- \zeta} \big)^{-1} ||_{(\alpha, \beta)} 
 \end{displaymath}
 then in view of the first line we get
 \begin{displaymath}
 || \big( \widetilde{z- \zeta} \big) ^{-1} ||_{(1, 0)} \leq  K_2(\alpha, \beta)^2 \frac{1}{ || \big( \widetilde{z- \zeta}  \big)  ||_{(1, 0 )}} =  K_2(\alpha, \beta)^2 \frac{1}{ || z- \zeta    ||_{(1, 0 )}}.
 \end{displaymath}
 Therefore the map $\big( \widetilde{z- \zeta} \big) ^{-1}$ is weakly singular in the ordinary sense.
 \qedhere
\begin{lemma} For every complex algebra of elliptic type and every $\varepsilon > 0$ 
\begin{displaymath}
\int\limits_{||z-\zeta||_{(1, 0)} = \varepsilon} \frac{d \tilde z}{\widetilde {z - \zeta}} =  2 \pi \hat{i}  = 2 \pi \frac{\beta + 2 i }{\sqrt{4 \alpha - \beta^2}} ,
\end{displaymath}
where $\tilde{z} = y - ix$ and $d \tilde z = dy - i dx$.
\end{lemma}
\proof We first show that for every $\varepsilon > 0$:
\begin{displaymath}
\int\limits_{||\widetilde{z-\zeta}||_{(\alpha, \beta)} = \varepsilon} \frac{d \tilde z}{\widetilde {z - \zeta}} =  2 \pi \hat{i}.
\end{displaymath}
We parametrize  $||\widetilde{z-\zeta}||_{(\alpha, \beta)} = \varepsilon$ as
\begin{displaymath}
\widetilde{z - \zeta} = \varepsilon \big( \cos(\theta) + \sin(\theta) \frac{\beta}{\sqrt{4 \alpha - \beta^2}} + i  \sin(\theta) \frac{2}{\sqrt{4 \alpha - \beta^2}} \big) = \varepsilon \big( \cos(\theta) + \hat{i} \sin(\theta) \big)
\end{displaymath}
which is a positive oriented curve if $\theta$ varies between $0$ and $2 \pi$. Hence by differentiating we get
\begin{displaymath}
d \tilde{z} = \varepsilon \hat{i} \big( \cos(\theta) + \hat{i} \sin(\theta) \big) d \theta,
\end{displaymath}
and thus
\begin{displaymath}
\int\limits_{||\widetilde{z-\zeta}||_{(\alpha, \beta)} = \varepsilon} \frac{d \tilde z}{\widetilde {z - \zeta}} = \int_{0}^{2 \pi} \frac{ \varepsilon \hat{i} \big( \cos(\theta) + \hat{i} \sin(\theta) \big) }{\varepsilon  \big( \cos(\theta) + \hat{i} \sin(\theta) \big)} d \theta =   2 \pi \hat{i} = 2 \pi \frac{\beta + 2 i }{\sqrt{4 \alpha - \beta^2}}.
\end{displaymath}
Let $\Omega$ be a bounded domain  with sufficiently smooth boundary $\partial \Omega$ and $g(z)$ a map valued in the complex algebra of elliptic type and that is continuously differentiable with respect to variables $x$ and $y$ in $\overline{\Omega}$. 
The Green-Gauss corresponding formulas are
\begin{displaymath}
\iint\limits_{\Omega} \partial_x g dx dy = \int\limits_{\partial \Omega} g dy,
\end{displaymath} 
\begin{displaymath}
\iint\limits_{\Omega} \partial_y g dx dy = -\int\limits_{\partial \Omega} g dx.
\end{displaymath} 
Multiplying the second line times $i$ and summing to the first line we get the complex Green-Gauss Integral Formula:
\begin{displaymath}
\iint\limits_{\Omega} \partial_{\bar z} g  \ dx dy = \frac{1}{2} \int\limits_{\partial \Omega } g d \tilde{z}.
\end{displaymath}
where $d \tilde{z} = dy - i dx $. \\
Let us show now that the circular curve $||z-\zeta||_{(1, 0)} =  \varepsilon $ lies in the interior of  the ellipse
\begin{displaymath}
 ||\widetilde{z-\zeta}||_{(\alpha, \beta)} =  \frac{\varepsilon +1}{K_{1}(\alpha, \beta)} .
 \end{displaymath}
From the comparison of the norms we have that for all $z$ and fixed $\zeta$
\begin{displaymath}
 ||\widetilde{z-\zeta}||_{(\alpha, \beta)} \leq \frac{1}{K_1(\alpha, \beta)} || \widetilde{z-\zeta}||_{(1, 0)}  =  \frac{1}{K_1(\alpha, \beta)} || z-\zeta||_{(1, 0)} .
\end{displaymath}
In particular, if $|| z-\zeta||_{(1, 0)}  = \varepsilon$ then
\begin{displaymath}
 ||\widetilde{z-\zeta}||_{(\alpha, \beta)} \leq \frac{\varepsilon}{K_1(\alpha, \beta)} < \frac{\varepsilon + 1}{K_1(\alpha, \beta)}.
\end{displaymath}
Now we set $\Omega_0$ as the region between the circular and ellipsoidal curves and $g(z)$ as the map $(\widetilde{z  - \zeta})^{-1}$. With the usual positive orientation of the Green-Gauss complex formula we get
\begin{displaymath}
0 = \frac{1}{2} \int\limits_{  ||\widetilde{z-\zeta}||_{(\alpha, \beta)} =  \frac{\varepsilon +1}{K_{1}(\alpha, \beta)}} \frac{ d \tilde{z}}{\widetilde{(z - \zeta)}} - \frac{1}{2} \int\limits_{  ||z-\zeta||_{(1, 0)} = \varepsilon} \frac{ d \tilde{z}}{\widetilde{(z - \zeta)}}
\end{displaymath}
since the map $(\widetilde{z - \zeta})^{-1}$ is holomorphic on $\overline{\Omega_0}$. Hence
\begin{displaymath}
\int\limits_{  ||z-\zeta||_{(1, 0)} = \varepsilon} \frac{ d \tilde{z}}{\widetilde{(z - \zeta)}}  = \int\limits_{  ||\widetilde{z-\zeta}||_{(\alpha, \beta)} =  \frac{\varepsilon +1}{K_{1}(\alpha, \beta)}} \frac{ d \tilde{z}}{\widetilde{(z - \zeta)}}.
\end{displaymath}
for every $\varepsilon > 0$. Setting
\begin{displaymath}
\varepsilon ' = \frac{\varepsilon +1}{K_{1}(\alpha, \beta)} > 0
\end{displaymath}
we have that 
\begin{displaymath}
 \int\limits_{  ||\widetilde{z-\zeta}||_{(\alpha, \beta)} = \varepsilon ' } \frac{ d \tilde{z}}{\widetilde{(z - \zeta)}} = 2 \pi \hat{i}.
\end{displaymath}
and it follows that for every $\epsilon>0$
\begin{displaymath}
\int\limits_{  ||z-\zeta||_{(1, 0)} = \varepsilon} \frac{ d \tilde{z}}{\widetilde{(z - \zeta)}} = 2 \pi \hat{i}.
\end{displaymath}
\qedhere
\begin{lemma}
For every  complex algebra of elliptic type and every continuous function $f$ on the region $|| \widetilde{z - \zeta}||_{(1, 0)} \leq \varepsilon$:
\begin{displaymath}
\lim_{\varepsilon \to 0} \int\limits_{|| z - \zeta||_{(1, 0)} = \varepsilon } \frac{f(z) - f(\zeta)}{\widetilde{z - \zeta}} d \tilde{z}= 0.
\end{displaymath}

\end{lemma}
\proof Using the triangle inequality with the norm depending on $\alpha$ and $\beta$ we get the following estimation
\begin{displaymath}
\big| \big| \int\limits_{|| z - \zeta||_{(1, 0)} = \varepsilon } \frac{f(z) - f(\zeta)}{\widetilde{z - \zeta}} d \tilde{z} \big| \big|_{(\alpha, \beta)} \leq \sup_{|| z - \zeta||_{(1, 0)} = \varepsilon} \frac{ || f(z) - f(\zeta)  ||_{(\alpha, \beta)}}{ || \widetilde{z - \zeta} ||_{(\alpha, \beta)}} \int\limits_{|| z - \zeta||_{(1, 0)}  = \varepsilon} \big| \big| d \tilde z  \big| \big|_{(\alpha, \beta)}
\end{displaymath}
where 
\begin{displaymath}
 \big| \big| d \tilde z  \big| \big|_{(\alpha, \beta)} = \sqrt{ \alpha dx^2 +  \beta dx dy + dy^2}.
\end{displaymath}
From the comparison of the norms we get
\begin{displaymath}
\sup_{|| z - \zeta||_{(1, 0)} = \varepsilon} \frac{1}{ || \widetilde{z - \zeta} ||_{(\alpha, \beta)}} \leq \sup_{|| z - \zeta||_{(1, 0)} = \varepsilon}  \frac{K_2(\alpha, \beta)}{ || \widetilde{z - \zeta} ||_{(1, 0)}} =   \frac{K_2(\alpha, \beta)}{ \varepsilon}
\end{displaymath}
and
\begin{displaymath}
  \sup_{|| z - \zeta ||_{(1, 0)}  = \varepsilon} || f(z) - f(\zeta)  ||_{(\alpha, \beta)} \leq  \frac{1}{K_{1}(\alpha, \beta)}   \sup_{|| z - \zeta ||_{(1, 0)}  = \varepsilon} || f(z) - f(\zeta)  ||_{(1, 0)}.
\end{displaymath}
Parametrizing the curve as $z = \varepsilon \big( \cos(\theta) + i \sin(\theta) \big) + \zeta$ and after differentiating with respect to $\theta$ we get
\begin{displaymath}
d \tilde{z} = \varepsilon \big( \cos(\theta) + i \sin(\theta) \big) d \theta.
\end{displaymath}  
Therefore
\begin{displaymath}
\int\limits_{|| z - \zeta||_{(1, 0)}  = \varepsilon} \big| \big| d \tilde z  \big| \big|_{(\alpha, \beta)}
\leq  2 \pi \varepsilon \sup_{0 \leq \theta < 2 \pi} || \cos(\theta) + i \sin(\theta) ||_{(\alpha, \beta)} 
\end{displaymath}
\begin{displaymath}
 \leq  \frac{2 \pi \varepsilon}{K_1(\alpha, \beta)} \sup_{0 \leq \theta < 2 \pi}  || \cos(\theta) + i \sin(\theta) ||_{(1, 0)} = \frac{2 \pi \varepsilon}{K_1(\alpha, \beta)}.
 \end{displaymath}
 In conclusion
 \begin{displaymath}
 \big| \big| \int\limits_{|| z - \zeta||_{(1, 0)} = \varepsilon } \frac{f(z) - f(\zeta)}{\widetilde{z - \zeta}} d \tilde{z} \big| \big|_{(\alpha, \beta)} \leq 2 \pi \frac{ K_2(\alpha, \beta)}{ K_1^2(\alpha, \beta)}  \sup_{|| z - \zeta ||_{(1, 0)}  = \varepsilon} || f(z) - f(\zeta)  ||_{(1, 0)}.
 \end{displaymath}
The continuity of the function $f$ implies that 
\begin{displaymath}
\sup_{||z - \zeta ||_{(1,0)} = \varepsilon} || f(z) - f(\zeta)  ||_{(1, 0)}
\end{displaymath}
goes to zero as $\varepsilon \to 0$. Therefore 
\begin{displaymath}
\lim_{\varepsilon \to 0} \big| \big| \int\limits_{|| z - \zeta||_{(1, 0)} = \varepsilon } \frac{f(z) - f(\zeta)}{\widetilde{z - \zeta}} d \tilde{z} \big| \big|_{(\alpha, \beta)}  = 0
\end{displaymath}
and hence
\begin{displaymath}
\lim_{\varepsilon \to 0} \int\limits_{|| z - \zeta||_{(1, 0)} = \varepsilon } \frac{f(z) - f(\zeta)}{\widetilde{z - \zeta}} d \tilde{z} = 0.
\end{displaymath}
\qedhere
\begin{lemma}
Let $g(z)$ a continuous  function with respect to variables $x$ and $y$ in $\overline{\Omega}$, a simply connected domain with finite measure but not
necessarily bounded. Suppose that
\begin{displaymath}
\sup_{\Omega} || g(z) ||_{(1, 0)}  < \infty. 
\end{displaymath}
Then the integral over $\Omega$ 
\begin{displaymath}
 \iint\limits_{\Omega} \frac{ g(z)}{\widetilde{z - \zeta}} dx dy
\end{displaymath}
exists .
\end{lemma}
\proof
 In order to prove this we consider first the following estimation:
\begin{displaymath}
|| \iint\limits_{\Omega} \frac{g(z)}{\widetilde{z - \zeta}} dx dy ||_{(\alpha, \beta)} \leq \sup_{\Omega} || g(z) ||_{(\alpha, \beta)} \iint\limits_{\Omega}   ||  (\widetilde{z - \zeta})^{-1}  ||_{(\alpha, \beta)} dx dy 
\end{displaymath} 
\begin{displaymath}
\leq \sup_{\Omega} || g(z) ||_{(1, 0)} \frac{1}{K_{1}^2(\alpha, \beta)} \iint\limits_{\Omega}   ||  (z - \zeta)^{-1}  ||_{(1, 0)} dx dy 
\end{displaymath}
then, in view of the second lemma we conclude that
\begin{displaymath}
|| \iint\limits_{\Omega} \frac{g(z)}{\widetilde{z - \zeta}} dx dy ||_{(\alpha, \beta)} \leq 
\sup_{ \Omega} || g(z) ||_{(1, 0)} \frac{K_{2}^2(\alpha, \beta)}{K_{1}^2(\alpha, \beta)} \iint\limits_{ \Omega} 
\frac{1}{ || (z - \zeta)  ||_{(1, 0)}} dx dy .
\end{displaymath}

The Schmidt inequality is given in this case by
\begin{displaymath}
\iint\limits_{ \Omega} \frac{1}{ || (z - \zeta)  ||_{(1, 0)}} dx dy \leq 2 \pi \big( \frac{m \Omega}{\pi})^{\frac{1}{2}},
\end{displaymath}
where $m \Omega$ denotes the measure of $\Omega$, which was assumed to be finite. Since
\begin{displaymath}
|| \iint\limits_{\Omega} \frac{g(z)}{\widetilde{z - \zeta}} dx dy ||_{(1, 0)} \leq \frac{1}{K_{1}(\alpha, \beta )} || \iint\limits_{\Omega} \frac{g(z)}{\widetilde{z - \zeta}} dx dy ||_{(\alpha, \beta)} 
\end{displaymath}
we have the desired result.
\qedhere
\section{Proof of the Cauchy-Pompeiu Formula}
Let $\Omega$ a simply connected and bounded domain with sufficiently smooth boundary.

We now introduce the domain
\begin{displaymath}
\Omega_{\varepsilon} = \Omega \setminus \overline{U_{\varepsilon} (\zeta)},
\end{displaymath}
where
\begin{displaymath}
U_{\varepsilon}(\zeta) =  \{  \ z  \ \ ; \ \  || z - \zeta ||_{(1, 0)} < \varepsilon \ \}.
\end{displaymath}

Setting $g(z)$ as
\begin{displaymath}
g(z) = \frac{f(z)}{\widetilde{z - \zeta}},
\end{displaymath}
where $f(z)$ is a continuously differentiable function with respect to variables $x$ and $y$ on $\overline{\Omega}$.

Since $(\widetilde{z-\zeta})^{-1}$ is continuosly differentiable except at $\zeta$, then $g$ 
is also continuosly differentiable except at $\zeta$.
Applying the Cauchy-Riemann operator to $g$ yields
\begin{displaymath}
\partial_{\bar z} g = \frac{\partial_{\bar z}f}{\widetilde{z - \zeta}}
\end{displaymath}
as the Cauchy-Riemann operator satisfies the product rule and the map $(\widetilde{z - \zeta})^{-1}$ is holomorphic.
The complex Green-Gauss integral formula on $\Omega_{\varepsilon}$ has then the form:
\begin{displaymath}
\iint\limits_{ \Omega_{\varepsilon}} \frac{\partial_{\bar z} f(z)}{\widetilde{z - \zeta}} dx dy = \frac{1}{2 } \int\limits_{\partial \Omega} \frac{f(z)}{\widetilde{z - \zeta}} d \tilde z - \frac{1}{2} \int\limits_{||z - z_0||_{(1,0)} = \varepsilon} \frac{f(z)}{\widetilde{z - \zeta}} d \tilde z,
\end{displaymath}
which we rewrite as
\begin{displaymath}
\iint\limits_{ \Omega_{\varepsilon}} \frac{\partial_{\bar z} f(z)}{\widetilde{z - \zeta}} dx dy = \frac{1}{2 } \int\limits_{\partial \Omega} \frac{f(z)}{\widetilde{z - \zeta}} d \tilde z - \frac{1}{2} \int\limits_{||z - z_0||_{(1, 0)} = \varepsilon} \frac{f(z) - f(\zeta) + f(\zeta)}{\widetilde{z - \zeta}} d \tilde z.
\end{displaymath}

Taking the limit as $\varepsilon \to 0$  on the left hand side we have:
\begin{displaymath}
\lim_{\varepsilon \to 0} \iint\limits_{ \Omega_{\varepsilon}} \frac{\partial_{\bar z} f(z)}{\widetilde{z - \zeta}} dx dy =  \iint\limits_{ \Omega} \frac{\partial_{\bar z} f(z)}{\widetilde{z - \zeta}} dx dy.
\end{displaymath}

In view of the precedent lemmas, and taking the limit as $\varepsilon \to 0$ for the right hand side also and solving in terms of $f(\zeta)$ the following theorem is proved:
\begin{theorem}
Suppose $f$ is continuosly differentiable on $\overline{\Omega}$. Then
\begin{displaymath}
f(\zeta) = \frac{1}{ 2 \pi \hat{i} }  \int\limits_{\partial \Omega} \frac{f(z)}{\widetilde{z - \zeta}} d \tilde z - \frac{1}{ \pi \hat{i}}  \iint\limits_{ \Omega} \frac{\partial_{\bar z} f(z)}{\widetilde{z - \zeta}} dx dy .
\end{displaymath}
\end{theorem}
As an consequence, we obtain the Cauchy integral formula for holomorphic functions:
\begin{corollary}
Suppose f is holomorphic on $\overline{\Omega}$ . Then
\begin{displaymath}
f(\zeta) = \frac{1}{ 2 \pi \hat{i} }  \int\limits_{\partial \Omega} \frac{f(z)}{\widetilde{z - \zeta}} d \tilde z.
\end{displaymath}
\end{corollary}
\section{Conclusions}
In this paper we have proved a generalization of the Cauchy-Pompeiu representation formula 
and also a Cauchy integral formula for holomorphic functions in parameter-depending complex 
numbers of elliptic type. This result suggests that a generalization of the ordinary theory of 
holomorphic functions can be developed. In particular the Cauchy representation formula 
allows to estimate the partial derivatives of holomorphic functions inside a ball which 
in turn is useful to generalize the study of the solvability of the initial value problems \cite{tut}.

\end{document}